\documentclass[11pt,lettersize,reqno,fullpage]{amsart}
\usepackage{amsmath}

\usepackage{float}

\usepackage{amsmath,amsfonts,amssymb}   
\usepackage{float}
\usepackage{amsmath,amssymb,amsfonts,mathrsfs,verbatim,enumitem,pstricks,amsthm, graphicx}
\usepackage[all]{xy}
\usepackage{centernot, xr, accents}
\usepackage{hyperref}

\usepackage{footnote}

\usepackage{hyperref}
\hypersetup{
	colorlinks,
	citecolor=red,
	filecolor=black,
	linkcolor=blue,
	urlcolor=black
}

\providecommand \hf{\hspace*{0.5cm}}
 
\newtheorem{thm}{Theorem}[section]

\newtheorem{conj}[thm]{Conjecture}
\newtheorem{lmm}[thm]{Lemma}

\newtheorem{defn}[thm]{Definition}
\newtheorem{prp}[thm]{Proposition}

\usepackage{amsmath,stackengine,scalerel}

\makeatletter
\renewcommand{\tocsection}[3]{%
  \indentlabel{\@ifnotempty{#2}{\bfseries\ignorespaces#1 #2\quad}}\bfseries#3}
\renewcommand{\tocsubsection}[3]{%
  \indentlabel{\@ifnotempty{#2}{\ignorespaces#1 #2\quad}}#3}

\newcommand\@dotsep{4.5}
\def\@tocline#1#2#3#4#5#6#7{\relax
  \ifnum #1>\c@tocdepth 
  \else
    \par \addpenalty\@secpenalty\addvspace{#2}%
    \begingroup \hyphenpenalty\@M
    \@ifempty{#4}{%
      \@tempdima\csname r@tocindent\number#1\endcsname\relax
    }{%
      \@tempdima#4\relax
    }%
    \parindent\z@ \leftskip#3\relax \advance\leftskip\@tempdima\relax
    \rightskip\@pnumwidth plus1em \parfillskip-\@pnumwidth
    #5\leavevmode\hskip-\@tempdima{#6}\nobreak
    \leaders\hbox{$\m@th\mkern \@dotsep mu\hbox{.}\mkern \@dotsep mu$}\hfill
    \nobreak
    \hbox to\@pnumwidth{\@tocpagenum{\ifnum#1=1\bfseries\fi#7}}\par
    \nobreak
    \endgroup
  \fi}
\AtBeginDocument{%
\expandafter\renewcommand\csname r@tocindent0\endcsname{0pt}
}
\def\l@subsection{\@tocline{2}{0pt}{2.5pc}{5pc}{}}
\makeatother

\begin{document}

\title[Weyl's Law for compact rank one symmetric spaces]{Weyl's Law for compact rank one symmetric spaces}


\author[Sai Sriharsha Indukuri]{Sai Sriharsha Indukuri}
\address{Department of Mathematics, Washington University in St.Louis, St.Louis 63105, Missouri, USA.}
\email{i.saisriharsha@wustl.edu}

\author[Ritwik Mukherjee]{Ritwik Mukherjee}
\address{School of Mathematical Sciences, National Institute of Science Education and Research, Bhubaneswar, An OCC of Homi Bhabha National Institute,
Khurda 752050, Odisha, India.}
\email{ritwikm@niser.ac.in}



\date{}

\maketitle

\begin{abstract}
Weyl's law is a fundamental result governing the asymptotic behaviour of 
the eigenvalues of the Laplacian.   
It states that 
for a compact $d$ dimensional manifold $M$ 
(without boundary), the eigenvalue counting function $N_M(\lambda)$
has the following asymptotic growth   
\[ N_M(\lambda) = \textnormal{Vol}(M)\left(\frac{\omega_d}{(2\pi)^d}\right) \lambda^d + O(\lambda^{d-1}). \]   
\hf \hf A natural question is: when is the error term sharp and when can it be improved? 
It has been known for a long time that the error term is sharp for the round sphere (since $1968$). 
In contrast, it has only recently been shown (in $2019$)  
by Iosevich and Wyman that for the product of spheres, the error term can be polynomially improved. 
They conjecture that a polynomial improvement should be true for products in general. \\ 
\hf \hf In this paper we extend both these results 
to Compact Rank One Symmetric Spaces (CROSSes). We show that for CROSSes, the error term is sharp. 
Furthermore, we show that for a product of CROSSes, 
the error term can be polynomially improved. This gives further evidence to the conjecture made by  
Iosevich and Wyman. 
\end{abstract}

\tableofcontents

\section{Introduction}
Let $(M,g)$ be a compact $d$-dimensional Riemannian manifold without boundary 
and let $\Delta_M$ be the  
Laplace-Beltrami operator (acting on functions). Associated to this operator,  
we define the eigenvalue counting function
\begin{align*}
N_{M}(\lambda)&:= \#\{ \textnormal{eigenvalues of $\Delta_M$} \leq \lambda^2 \}.
\end{align*}   
In other words, $N_{M}(\lambda)$ is the number of eigenvalues of $\Delta_M$ less than or equal to $\lambda^2$ 
(counted with multiplicities). It is a fundamental result that the asymptotic growth of $N_{M}(\lambda)$ 
is governed by Weyl's law: 
\begin{align} 
\lim_{\lambda\longrightarrow \infty}\frac{N_{M}(\lambda)}{\lambda^d} & = \textnormal{Vol}(M)\left(\frac{\omega_d}{(2\pi)^d}\right), 
\label{Weyl_leading_Term}
\end{align}
where 
$\omega_d$ is the volume of the unit ball in $\mathbb{R}^d$. A proof of \eqref{Weyl_leading_Term} 
can be found in \cite[page 95, corollary 2.43]{Getz_Ber}; it uses the the asymptotic expansion of the heat kernel. \\ 
\hf \hf As with any question on asymptotics, it is natural to ask what can be said about the error term. 
In $1920$, Richard Courant showed that the error term is of the order of $O(\lambda^{d-1} \log(\lambda))$ 
(\cite{R_Courant}). Since then, an extensive amount of work has been done by mathematicians to 
study the error term. \\ 
\hf \hf It was shown by Lars H{\"o}rmander in $1968$ (\cite{LH}) that the error term is  
$O(\lambda^{d-1})$. 
He further showed in that paper that the $O(\lambda^{d-1})$ 
error term is sharp for the round sphere $S^d$ 
(the special case of the three dimensional sphere $S^3$ was worked out earlier by Avacumovi\v{c} in \cite{Ava}). \\ 
\hf \hf Although the example of the sphere shows that the error, in general, can not 
be made better than $O(\lambda^{d-1})$, it is still 
natural to ask under what special cases the error can be improved. In \cite{Dust} and \cite{Ivrii}, Duistermaat and Guillemin 
and Ivrii show that if the measure of all periodic geodesic billiards is zero, then the error term can be improved 
(see \cite[page 595, theorem 1.3]{Weyl_Wyman_Iosevich} for the exact statement). \\ 
\hf \hf In \cite{Weyl_Wyman_Iosevich}, Iosevich and Wyman show 
\cite[page 596, theorem 1.4]{Weyl_Wyman_Iosevich}
that for a product of 
at least two compact 
Riemannian manifolds without boundary, the 
hypothesis of 
Duistermaat and Guillemin 
and Ivrii's Theorem holds; hence the error term can be improved. But that is not all, they show that 
for the product of spheres, the 
error term 
can be \textit{polynomially} improved. More precisely, they show \cite[page 596, theorem 1.4]{Weyl_Wyman_Iosevich} 
that if $M$ is the product of $n$ spheres, then 
\begin{align}
N_M(\lambda) & = \textnormal{Vol}(M)\left(\frac{\omega_d}{(2\pi)^d}\right)  \lambda^d + O(\lambda^{d-1-\frac{n-1}{n+1}}).
\end{align}
Based on this theorem, Iosevich and Wyman make the following conjecture. 
\begin{conj}\label{Conj_IW} 
\cite[page 597]{Weyl_Wyman_Iosevich} 
Let $M:= M_1 \times \ldots \times M_n$ be the product of 
at least two compact Riemannian manifolds without boundary. Then there exists a $\delta>0$ 
such that  
\begin{align}
N_M(\lambda) & = \textnormal{Vol}(M)\left(\frac{\omega_d}{(2\pi)^d}\right)  \lambda^d + O(\lambda^{d-1-\delta}).
\end{align}
\end{conj}
In this paper, we give further evidence to support the above conjecture. \\ 
\hf \hf Before stating the main result of this paper, let us 
recall from \cite[page 76]{AB} 
that a Compact Rank One Symmetric Spaces (CROSS) 
is either the sphere $S^d$, the 
real projective space $\mathbb{RP}^d$, 
the complex projective space $\mathbb{CP}^d$, the quaternionic projective space $\mathbb{HP}^d$ 
or the two dimensional octonionic projective space $\mathbb{OP}^2$. 
Note that aside from $S^1$ and $\mathbb{RP}^1$ (which is again $S^1$ with a  different radius), all the other CROSSes 
have real dimension at least two. 
The main result of this paper is as follows. 
\begin{thm}
\label{MT_prod}
Let $M:= M_1 \times \ldots \times M_n$ be the product of 
at least two compact Riemannian manifolds without boundary. 
Suppose each $M_i$ is one of the CROSSes of real dimension at least two. 
Then 
\begin{align}
N_M(\lambda) & = \textnormal{Vol}(M)\left(\frac{\omega_d}{(2\pi)^d}\right)  \lambda^d + O(\lambda^{d-1-\frac{n-1}{n+1}}).
\end{align}
\end{thm}
Note that theorem \ref{MT_prod} supports Conjecture \ref{Conj_IW}, with $\delta:= \frac{n-1}{n+1}$. \\ 
\hf \hf The second result of our paper is as follows.   
\begin{thm}
\label{MR2}
Let $M$ be one the 
$CROSSes$.
Then the error term in Weyl's law can not be improved. 
\end{thm}


\section{Structure of the paper} 
In order to prove theorem \ref{MT_prod}, we use a fundamental 
lattice counting result proved by Iosevich and Wyman 
\cite[page 597, theorem 1.8]{Weyl_Wyman_Iosevich}. 
The main observation of this paper is that their lattice counting result   
applies to a much more general setting than spheres; in particular it applies to all CROSSes. \\   
\hf \hf In section \ref{setup_not}, we describe the general setup in which we can prove  
theorem \ref{MT_prod}. The approach will work if the eigenvalues and multiplicities are explicitly given  
and a certain condition on the multiplicity function can be verified. \\ 
\hf \hf In section \ref{Error_can_be_improved} we show that 
theorem \ref{MT_prod} is true for product of all manifolds whose multiplicity function satisfies the symmetry condition. 
In the proof, we make extensive use of Iosevich and Wyman's lattice counting result. 
In section \ref{CROSS_is_W}, 
we verify that the symmetry condition is indeed satisfied for all the CROSSes. \\ 
\hf \hf Finally in section \ref{Error_is_sharp}, we show that the error is sharp for a single CROSS.

\section{Setup and notation}
\label{setup_not}
Our set up is as follows. 
Let $(M, g)$ be a compact Riemannian manifold without boundary, with dimension $d \geq 2 $.
We consider the eigenvalues of the Laplacian acting on functions. Assume that the (square root of) the 
eigenvalues are given by 
\begin{align}
\Lambda_M&:= \{ \sqrt{\alpha_M(k)}: k =0, 1, 2, \ldots \}, \qquad \textnormal{where} \qquad 
\alpha_{M}(k):= A_M k^2 + B_M k + C_M. \label{Lambda_condt}
\end{align} 
Here $A_M, B_M$ and $C_M$ are constants depending on $M$. 
We will assume that $A_M$ is not zero, i.e.  
\begin{align}
A_M&>0. \label{A_condt}
\end{align}
Assume that the multiplicity of these eigenvalues are given by a degree $d-1$ polynomial
$\mathsf{R}_{M}(k)$.
Define 
\begin{align}
\mathsf{Q}_{M}(t)&:= \mathsf{R}_{M}\Big(t-\frac{B_M}{2 A_M}\Big). \label{QM_defn}
\end{align} 
We now state the crucial assumption that the multiplicity function $\mathsf{Q}_M$ should satisfy: 
we assume that 
\begin{align}
\mathsf{Q}_{M}(t) & = C t^{d-1} + \varphi(t), \qquad \textnormal{if} \qquad d>2 \qquad \textnormal{and} \nonumber \\ 
\mathsf{Q}_{M}(t) & = C t, \qquad \textnormal{if} \qquad d=2,\label{Q_condt}
\end{align}
where $\varphi(t)$ is a polynomial of degree at most $d-3$ and 
$C$ is some constant. 

\begin{defn} 
Let $(M, g)$ be a compact Riemannian manifold without boundary, with dimension $ d \geq 2$. 
If conditions \eqref{Lambda_condt}, \eqref{A_condt} and \eqref{Q_condt} are satisfied, then 
we will say that 
$(M, g)$ 
is 
a $\mathsf{W}$-manifold. 
\end{defn}
We have used the letter $\mathsf{W}$ in honour of Weyl. \\
\hf \hf The condition of being a $\mathsf{W}$-manifold 
will play a crucial role in the proof of theorem \ref{MT_prod}.  
It will also play a role in the proof of theorem \ref{MR2}.
The main results of this paper can be rephrased by the following three theorems.   

\begin{thm}
\label{main_thm_W}
Let $(M_1, g_1), (M_2, g_2), \ldots, (M_n, g_n)$ be 
$\mathsf{W}$-manifolds of dimensions $d_1, d_2, \ldots, d_n$. 
Let 
\begin{align}
M&:= M_1 \times M_2 \times \ldots \times M_n, \qquad \textnormal{and} \qquad d:= d_1 + d_2 \ldots + d_n. \nonumber
\end{align} 
If $n>1$, then 
\begin{align}
N_M(\lambda) & = \textnormal{Vol}(M)\left(\frac{\omega_d}{(2\pi)^d}\right) \lambda^d + O\left(\lambda^{d-1-\frac{n-1}{n+1}}\right). \nonumber
\end{align}
\end{thm} 

\begin{thm} 
\label{main_thm_W2_sharp}
Let $(M, g)$ be a  
$\mathsf{W}$-manifold. Then the error term in Weyl's law can not be improved.
\end{thm} 
 
\begin{thm}
\label{CROSS_is_W_ag}
Compact Rank One Symmetric Spaces 
of real dimension $d\geq 2$ 
are $\mathsf{W}$-manifolds.   
\end{thm}
We note that theorems 
\ref{main_thm_W}, \ref{main_thm_W2_sharp} and \ref{CROSS_is_W_ag} 
immediately imply theorems \ref{MT_prod} and \ref{MR2}. 


\section{A lattice counting result by Iosevich and Wyman} 
In this section, we recall the lattice counting theorem proved by Iosevich and Wyman in \cite{Weyl_Wyman_Iosevich}.  
\begin{thm}
\label{IW_Main_Result}
\cite[page 597, theorem 1.8]{Weyl_Wyman_Iosevich}
Fix $y \in \mathbb{R}^{n} $ and let 
$(d_{1},d_{2},\dots,d_{k},1,\dots,1) \in \mathbb{N}^{n} $ be a multiindex with $d_{i}\geq 2 $ 
for $i=1,\dots,k $. Let $p$ be a norm on $\mathbb{R}^{n} $ such that $p^{-1}(1) $ is a smooth hypersurface in $\mathbb{R}^{n} $ with nonvanishing Gaussian curvature. Define $d:= d_1 + \ldots + d_{k}+ n-k$.
For a weight function \[F(x)=\prod_{i=1}^{k} \chi_{[0,\infty)}(x_{i})\ x_{i}^{d_{i}-1}, \]
we have \[\sum_{{m \in (\mathbb{Z}^{n}+y)}\atop {p(m)\leq \lambda}} F(m)=\lambda^{ d} \int_{p^{-1}([0,1])} F(x) \ dx+E(\lambda),\] 
where 
$E(\lambda) $ satisfies the bounds \[E(\lambda)=O\left(\lambda^{ d-1-\frac{n-1}{n+1}} \right) \]
uniformly in $y$.
\end{thm}
Using this lattice counting theorem, we are going to prove the following result. 
\begin{thm}
\label{IW_Lattice}
Let $d_1, d_2, \ldots, d_n \in \mathbb{Z}^{\geq 2}$ and $n\geq 2$. 
Let $y_1, \ldots, y_n$ be fixed $n$ non negative real numbers and 
$P_1, P_2, \ldots P_n$ be fixed $n$ positive real numbers. 
Then
\begin{align}
\sum_{\substack{m_1-y_1, \ldots, m_n-y_n \in \mathbb{Z}^{\geq 0}\\ P_1 m_1^2 + \ldots + P_n m_n^2 \leq \lambda^2}} 
m_1^{d_1-1}\ldots m_n^{d_n-1}&= C \lambda^{d} 
+ O\left(\lambda^{d-1-\frac{n-1}{n+1}}\right),
\end{align}
where $d:= d_1 + d_2 \ldots + d_n$ and $C$ is a constant. 
\end{thm}
We first prove an auxiliary result 
that combined with Theorem \ref{IW_Main_Result} will immediately imply 
Theorem \ref{IW_Lattice}. First we make the following notation: define 
the norm $p$ on $\mathbb{R}^{n} $ by 
\[p(m):= \left( P_{1}m_{1}^{2}+\dots+P_{n}m_{n}^{2}\right)^{\frac{1}{2}},\ m=(m_{1},\dots,m_{n})\in \mathbb{R}^{n}. \] 
The crucial fact that will eventually be used (in order to apply Theorem \ref{IW_Main_Result}) 
is that 
the ellipsoid $p^{-1}(1) $
has non vanishing Gaussian curvature.
We will denote $B_p$ to be the unit ball in $\mathbb{R}^n$ with respect to the norm $p$, i.e.   
\begin{align*}
B_p &:= \{\xi \in \mathbb{R}^n : p(\xi)^2 \leq 1\}. 
\end{align*} 

\begin{lmm}
\label{Int_Lmm} 
With notation as above, define 
\begin{align*}
N(\lambda) &:= \sum_{{m \in (\mathbb{Z}_{\geq 0}^{n} + y)}
\atop{p(m)^2 \leq \lambda^2}}m_1^{d_1-1} \ldots m_n^{d_n-1} 
\qquad \textnormal{and} \qquad 
N_2(\lambda) := \sum_{{m \in (\mathbb{Z}^n + y) 
\cap \mathbb{R}_{\geq 0}^{n}}\atop{p(m)^2 \leq \lambda^2}}m_1^{d_1-1}\ldots m_n^{d_n-1}.
\end{align*} 
Then, 
\begin{align}
N_2(\lambda)&= N(\lambda) + O(\lambda^{d-2}). \label{N2_eq_N_plus_O_lm_d_minus_2}
\end{align}
\end{lmm}
\noindent\textbf{Proof:} 
Note that we can write $N(\lambda)$ as 
\begin{align*}
N(\lambda) & = \sum_{m \in \mathbb{Z}^{n}_{\geq 0} + y} \chi_{B_p} (m/\lambda)\ m_1^{d_1-1} \ldots m_n^{d_n-1}, 
\end{align*}
where $\chi$ denotes the characteristic function. 
Let us define the following sets for $j \in \{1, 2, \ldots, n\}$
\begin{align*} 
S&:= \big(\mathbb{Z}^n+y\big)\cap \mathbb{R}_{\geq 0}^n, \qquad
H_j:= \{x \in \mathbb{R}^n_{\geq 0} : x_j \in [0, y_j + 2)\}, \\ 
L_j&:= \big(\mathbb{Z}^n+y\big)\cap H_j  
\qquad \textnormal{and} \qquad A_j:= \{m \in \mathbb{Z}_{\geq 0}^n+y : m_j-y_j \in \{0, 1\}\}.
\end{align*}
Using the fact that $y_j \geq 0$ for all $j$, 
the following two set theoretic identities can easily be verified 
\begin{align} 
S&= \Big(\mathbb{Z}^{n}_{\geq 2} + y\Big) \bigsqcup \bigcup_{j=1}^{n} L_j 
\qquad \textnormal{and} \label{S_Tj_du} \\ 
\big(\mathbb{Z}_{\geq 0}^n+y\big) & = \big(\mathbb{Z}_{\geq 2}^n+y\big) \bigsqcup \bigcup_{j=1}^{n}A_j, 
\label{Aj_du} 
\end{align} 
where $\sqcup$ denotes disjoint union. Note that 
$y_j \geq 0$ is used only to prove \eqref{S_Tj_du} 
(for \eqref{Aj_du} that fact is not relevant). 
Let us now rewrite our given functions using the notation we have introduced. 
We note that 
\begin{align} 
N(\lambda) &= \sum_{m \in \mathbb{Z}^n_{\geq 0} + y} \chi_{B_p} (m/\lambda)\ m_1^{d_1-1} \ldots m_n^{d_n-1} 
\qquad \textnormal{and} \qquad \label{N_redef}\\  
N_1(\lambda)& = \sum_{m \in \mathbb{Z}^n_{\geq 2} + y} \chi_{B_p} (m/\lambda)\ m_1^{d_1-1} \ldots m_n^{d_n-1} 
\label{N1_redef}
\end{align} 
Note that since $y_i\geq 0$ for all $i \in \{1, 2, \ldots, n\}$, we conclude that the 
summand on the right hand side of equation \eqref{N_redef} is non negative for all $m$ in the 
domain of the summation. Furthermore, $N_1$ is defined by summing up the same function over a smaller 
set. Hence, using \eqref{Aj_du}, we conclude that 
\begin{align} 
N_1(\lambda) \leq N(\lambda) & \leq N_1(\lambda) + \sum_{j=1}^{n} \sum_{m\in A_j} \chi_{B_p}(m/\lambda) 
m_1^{d_1-1} \ldots m_n^{d_n-1}. \label{N_leq_N1} 
\end{align}
Note that $y_i\geq 0$ is crucial for this argument since that ensures that the summand is non negative. 
Next, let us define for all $j \in \{1, \ldots, n\}$ 
and for all $\sigma \in \{0, 1\}$  
\begin{align} 
\mathsf{R}_{j\sigma}&:= 
\sum_{{m \in \mathbb{Z}^{n}_{\geq 0} + y}\atop {m_j = y_j+\sigma}}\chi_{B_p}(m/\lambda) 
\ m_1^{d_1-1} m_2^{d_2-1} \ldots m_n^{d_n-1}. \label{R_i_sigma}
\end{align}
We now note that the second term on the right hand side of  
\eqref{N_leq_N1} is precisely equal to 
\begin{align*}
\sum_{j=1}^{n} \sum_{\sigma=0}^{1} \mathsf{R}_{j\sigma}
\end{align*}
Hence, equation \eqref{N_leq_N1} can be rewritten as 
\begin{align}
0 \leq N(\lambda) - N_1(\lambda) &\leq \sum_{j=1}^{n} \sum_{\sigma=0}^{1} \mathsf{R}_{j\sigma} \nonumber \\ 
\implies \big|N(\lambda) - N_1(\lambda)\big| &\leq \sum_{j=1}^{n} \sum_{\sigma=0}^{1} \mathsf{R}_{j\sigma}. 
\label{N_minus_N1}
\end{align}
Next, let us study the function $N_2$. We note that 
\begin{align} 
N_2(\lambda)& = \sum_{m\in S} \chi_{B_p}(m/\lambda) m_1^{d_1-1} \ldots m_n^{d_n-1}.
\end{align} 
Using equation \eqref{S_Tj_du} and the non-negativity of the summand in the definition of $N_2$, 
we conclude that 
\begin{align}
0 \leq N_2(\lambda) - N_1(\lambda) &\leq \sum_{j=1}^{n} \sum_{m\in L_j}   
\chi_{B_p}(m/\lambda) m_1^{d_1-1} m_2^{d_2-1} \ldots m_n^{d_n-1} \nonumber \\ 
\implies \big|N_2(\lambda) - N_1(\lambda)\big| &\leq 
\sum_{j=1}^{n} \sum_{m\in L_j}   
\chi_{B_p}(m/\lambda) m_1^{d_1-1} m_2^{d_2-1} \ldots m_n^{d_n-1}.
\label{N2_minus_N1}
\end{align}
Hence, to prove the Lemma 
(namely equation \eqref{N2_eq_N_plus_O_lm_d_minus_2})
it suffices to show that   
it suffices to show that  
\begin{align}
&\mathsf{R}_{j\sigma}  = O(\lambda^{d-2}) \qquad \forall ~~j \in \{1, \ldots, n\} 
\qquad \textnormal{and} \qquad \forall \sigma \in \{0, 1\} \qquad \textnormal{and} \label{Rj_sigma_ts}\\ 
& \sum_{m \in L_j} \chi_{B_p} (m/\lambda)\ m_1^{d_1-1} \ldots m_n^{d_n-1} = O(\lambda^{d-2}) 
\qquad  \forall ~~j \in \{1, \ldots, n\}. \label{Lj_ts}
\end{align}
To see why that is so, we note that if we do manage to prove the above two claims, then 
by \eqref{N_minus_N1}, we will conclude that $N(\lambda)-N_1(\lambda) \in O(\lambda^{d-2})$ 
while using \eqref{N_minus_N1}, we will conclude that   
$N_2(\lambda)-N_1(\lambda) \in O(\lambda^{d-2})$. Using the fact that if $f , g \in O(\lambda^k)$ 
then $f-g$ also belongs to $O(\lambda^k)$, we conclude that 
$N_2-N = (N_2-N_1)-(N-N_1)$ is in $O(\lambda^{d-2})$. Let us now prove equations 
\eqref{Rj_sigma_ts} and \eqref{Lj_ts}.   \\
\hf \hf We start by proving \eqref{Rj_sigma_ts}.
Without loss of generality, it 
is enough to prove equation \eqref{Rj_sigma_ts} for $j=1$. 
We show that now. 
In the subsequent calculations, we will make a standard abuse of notation often used by analysts: if the letter 
$C$ appears in two different lines, then they usually refer to two different constants. 
Keeping this in mind, we note that 
\begin{align} 
\mathsf{R}_{1\sigma} & = (y_1 + \sigma)^{d_1-1}
\sum_{{m_2 - y_2, \ldots ,m_n-y_n \in \mathbb{Z}_{\geq 0}}\atop {m_1 = y_1+\sigma}}\chi_{B_p}(m/\lambda) \ m_2^{d_2-1} \ldots m_n^{d_n-1}, \label{R_1_zeroth}\\ 
\implies \mathsf{R}_{1\sigma} 
& \leq C\sum_{{m_2 - y_2, \ldots ,m_n-y_n \in \mathbb{Z}_{\geq 0}}\atop {m_1 = y_1+\sigma}}\chi_{B_p}(m/\lambda) \ p(m)^{d-d_1-n+1}. \label{R_1_sigma_leq_C_t} 
\end{align}
Let us 
justify 
\eqref{R_1_sigma_leq_C_t}.  
Note that for any $j \in \{2, \ldots, n\}$, we have 
\begin{align*} 
m_j^{2(d_j-1)}&= \Big(\Big(\frac{1}{P_j}\Big)^{2(d_j-1)}\Big) 
\Big((P_j m_j)^{2}\Big)^{d_j-1}. 
\end{align*} 
Hence, 
\begin{align} 
m_j^{2(d_j-1)}&\leq \Big(\frac{1}{P_j}\Big)^{2(d_j-1)} \Big(P_1^2m_1^2+ \ldots + P_n m_n^2\Big)^{d_j-1}  
 \qquad \forall~~j \in \{2, \ldots, n\}. \label{mj_pp}
\end{align}
Multiplying the inequalities \eqref{mj_pp} from $j=2$ to $n$ gives us 
\begin{align} 
m_2^{2(d_2-1)}m_3^{2(d_3-1)} \ldots m_n^{2(d_n-1)} &\leq C   
\Big(P_1^2m_1^2+ \ldots + P_n m_n^2\Big)^{d-d_1-n+1}. 
\label{prod_mj_llkk} 
\end{align}
Equation \eqref{R_1_zeroth} combined with 
equation \eqref{prod_mj_llkk} gives us 
equation \eqref{R_1_sigma_leq_C_t} 
(we remind the reader that the $C$ in equations \eqref{prod_mj_llkk} and \eqref{R_1_sigma_leq_C_t} 
refer to two different constants). 
Next, we will use the following standard fact about the asymptotic behaviour about the 
number of lattice points inside a ball of radius $\lambda$. 
Define  
\begin{align} 
M(\lambda)&:=  
\sum_{{m_1-y_1, m_2 - y_2, \ldots ,m_n-y_n \in \mathbb{Z}}\atop {P_1 m_1^2+\ldots + P_n m_n^2 \leq \lambda^2}} 
\mathbf{1}. \label{M_ball_std}
\end{align}
Then 
\begin{align} 
M(\lambda)& \leq C \lambda^n. \label{M_lm} 
\end{align}
Although this is a standard fact, we have included a proof of 
equation \eqref{M_lm} in section \ref{ball_std}. 

Using equation \eqref{M_lm}, we conclude that 
\begin{align} 
\mathsf{R}_{1\sigma}&\leq C \lambda^{d-d_1-n+1}
\sum_{{m_2 - y_2, \ldots ,m_n-y_n \in \mathbb{Z}_{\geq 0}}\atop {m_1 = y_1+\sigma}}\chi_{B_p}(m/\lambda) 
\nonumber \\ 
& = C \lambda^{d-d_1-n+1} \sum_{\substack{m_2-y_2, m_3 - y_3, \ldots ,m_n-y_n \in \mathbb{Z} 
\nonumber \\ 
P_2 m_2^2+\ldots + P_n m_n^2 \leq \lambda^2-P_1m_1^2 \\ m_1 = y_1 + \sigma}}  \mathbf{1} 
\nonumber \\ 
& \leq C \lambda^{d-d_1-n+1} \sum_{\substack{m_2-y_2, m_3 - y_3, \ldots ,m_n-y_n \in \mathbb{Z} \\ 
P_2 m_2^2+\ldots + P_n m_n^2 \leq \lambda^2}}  \mathbf{1} \nonumber \\ 
& \leq C \lambda^{d-d_1-n+1} \lambda^{n-1} \label{lattice_int_n_minus_1} \\ 
& = C \lambda^{d-d_1} \nonumber \\ 
& = O(\lambda^{d-2}) \qquad \textnormal{(since $d_1\geq 2$).} \label{d1_geq_2_O}
\end{align}
Equation \eqref{lattice_int_n_minus_1} follows from equation 
\eqref{M_lm} (with $n$ replaced by $n-1$).  
Equation \eqref{d1_geq_2_O} follows from the fact that $d_1\geq 2$. 
This proves \eqref{Rj_sigma_ts}. It now remains to show 
\eqref{Lj_ts}, which we now prove. 
Again, it suffices to prove \eqref{Lj_ts} for $j=1$.  
Note that 
\begin{align} 
\sum_{m \in L_1} \chi_{B_p} (m/\lambda)\ m_1^{d_1-1} \ldots m_n^{d_n-1} 
& \leq (y_1 + 2)^{d_1-1} \sum_{m \in L_1}\chi_{B_p}(m/ \lambda) m_2^{d_2-1} \ldots  m_n^{d_n-1} \nonumber \\ 
& \leq C \ \lambda^{d-d_1-n+1} \sum_{m \in L_1}\chi_{B_p}(m/ \lambda) \label{eq_lm_1}\\ 
& = C \ \lambda^{d-d_1-n+1} \lambda^{n-1} \label{eq_lm_2} \\ 
& = C \ \lambda^{d-d_1} \nonumber \\ 
& = O(\lambda^{d-2}). \label{eq_lm_3}
\end{align}
Equation \eqref{eq_lm_1} follows from equation \eqref{prod_mj_llkk}. Equation 
\eqref{eq_lm_2} follows from \eqref{M_lm} (with $n$ replaced with $n-1$). 
Finally, equation \eqref{eq_lm_3} follows from the fact that $d_1 \geq 2$. 
This proves Lemma \ref{Int_Lmm}. \qed  \\

\noindent\textbf{Proof of Theorem \ref{IW_Lattice}:} From Lemma \ref{Int_Lmm}: 
\begin{equation}
N_2(\lambda) = N(\lambda) + O(\lambda^{d-2}) \label{k1}
\end{equation}
and Theorem \ref{IW_Main_Result} 
gives:
\begin{equation}
N_2(\lambda) = C\lambda^{d} + O(\lambda^{d-1-\frac{n-1}{n+1}}). \label{k2}
\end{equation}
Equations \eqref{k1} and \eqref{k2}
combine to give: \(N(\lambda) = C\lambda^{d}  + O(\lambda^{d-1-\frac{n-1}{n+1}})\) 
and that is precisely the statement of Theorem \ref{IW_Lattice}. \qed

\section{Error can be polynomially improved for products of $\mathsf{W}$-manifolds: proof}
\label{Error_can_be_improved}
Let us first prove an auxiliary result. 
In the subsequent statements, $C$ always refers to a constant.  
We are going to make the standard abuse of notation in where we refer to different constants by the same letter $C$. 
\begin{thm}
\label{Auxiliary_result}
Let $d_1, d_2, \ldots, d_n \in \mathbb{Z}^{\geq 2}$ and $n\geq 2$. 
let $\mathsf{Q}_{d_i-1}(t)$ be a degree $d_i-1$ polynomial in $t$, 
that does not have a $t^{d_i-2}$ term, i.e.  
\begin{align}
\mathsf{Q}_{d_i-1}(t) & = C_i t^{d_i-1} + \varphi_i(t), \qquad \textnormal{if} \qquad d_i>2 \qquad \textnormal{and} \nonumber \\ 
\mathsf{Q}_{d_i-1}(t) & = C_i t, \qquad \textnormal{if} \qquad d_i=2,\label{Q_condt_ag}
\end{align}
where $C_i$ is a constant and $\varphi_i(t)$ is a polynomial of degree at most $d_i-3$.  
Let $y_1, \ldots, y_n$ be $n$ fixed real numbers and 
$P_1, P_2, \ldots P_n$ be $n$ fixed positive real numbers. 
Given a real number $\lambda$, let $N(\lambda)$ be defined as   
\begin{align}
N(\lambda) & = \sum_{\substack{m_1-y_1, \ldots, m_n-y_n \in \mathbb{Z}^{\geq 0}\\ P_1 m_1^2 + \ldots + P_n m_n^2 \leq a \lambda^2 + b}}
\mathsf{Q}_{d_1-1}(m_1) \ldots \mathsf{Q}_{d_n-1}(m_n),
\end{align} 
where $a \in \mathbb{R}^{>0}$ and $b \in \mathbb{R}$. Then,  
\begin{align}
N(\lambda) & = C \lambda^{d} + O(\lambda^{d-1-\frac{n-1}{n+1}}),
\end{align}
where $d:= d_1 + d_2 \ldots + d_n$ and $C$ is a constant.
\end{thm} 

\noindent \textbf{Proof:} Let us unwind the summand of $N(\lambda)$. Let us denote the summation of the leading 
power as $T_1$. The sum of the remaining powers will be denoted by $T_2$. In other words, 
\begin{align}
T_1 &:= \sum_{\substack{m_1-y_1, \ldots, m_n-y_n \in \mathbb{Z}^{\geq 0}\\ P_1 m_1^2 + \ldots + P_n m_n^2 \leq a\lambda^2+b}} 
C m_1^{d_1-1}\ldots m_n^{d_n-1}.
\end{align} 
Note that the sum total of the powers of all the $m_i$ is $d-n$. \\ 
\hf \hf Next, we note that $T_2$ is a finite sum of the type 
\begin{align} 
\sum_{\substack{m_1-y_1, \ldots, m_n-y_n \in \mathbb{Z}^{\geq 0}\\ P_1 m_1^2 + \ldots + P_n m_n^2 \leq a\lambda^2+b}} 
C m_1^{e_1}\ldots m_n^{e_n},
\end{align}
where the sum of the total power $e_1+\ldots e_n$ is at most $d-n-2$ 
(due to the condition \ref{Q_condt_ag}). In particular, the sum total of the powers of $m_i$ in the 
summand of each of the terms of $T_2$ is at least $2$ less than the 
sum total of the powers  of the $m_i$ in the summand of $T_1$. \\ 
\hf \hf Using theorem \ref{IW_Lattice}, we conclude that 
\begin{align}
T_1& = C \Big(a \lambda^2 + b\Big)^{\frac{d}{2}} + O\left((a \lambda^2 + b)^{\frac{1}{2}\left(d-1-\frac{n-1}{n+1}\right)}\right). \label{T1_1}
\end{align}
Using binomial theorem, we conclude that 
\begin{align} 
\Big(a \lambda^2 + b\Big)^{\frac{d}{2}} & = a^{\frac{d}{2}} \lambda^{d} + O(\lambda^{d-2}). \label{T1_2}
\end{align}
Combining \eqref{T1_1} and \eqref{T1_2}, we conclude that 
\begin{align}
T_1 & = C \lambda^{d} + O(\lambda^{d-1-\frac{n-1}{n+1}}). \label{T1_3}
\end{align}
Similarly, we again apply theorem \ref{IW_Lattice} and the 
binomial theorem to each of the terms of $T_2$. Doing that, we 
conclude that 
$T_2$ will be given by the right hand side of \eqref{T1_3} with $d$ replaced by $d-2$, i.e. 
\begin{align} 
T_2 & = C \lambda^{d-2} + O(\lambda^{d-3-\frac{n-1}{n+1}}). \label{T2_initial}
\end{align}
However, \eqref{T2_initial} implies that 
\begin{align} 
T_2 & = O(\lambda^{d-2}). \label{T2_final}
\end{align}
Since $N(\lambda) = T_1+T_2$, we conclude using \eqref{T1_3} and \eqref{T2_final}, that 
\begin{align}
N(\lambda) & = C \lambda^{d} + O(\lambda^{d-1-\frac{n-1}{n+1}}),
\end{align}
thereby proving theorem \ref{Auxiliary_result}. \qed\\  
\hf \hf We are now ready to prove theorem \ref{main_thm_W}. \\ 

\noindent \textbf{Proof of theorem \ref{main_thm_W}:} Let $N_M(\lambda)$ denote the eigenvalue counting function on $M$. 
Recall that $N_M(\lambda) $ can be expressed as 
\begin{align} 
N_M(\lambda) & = \sum_{\substack{k_1, \ldots, k_n \in \mathbb{Z}^{\geq 0}\\ 
 \alpha_{M_1}(k_1) + \ldots + \alpha_{M_n}(k_n) \leq \lambda^2}} \mathsf{R}_{M_1}(k_1) \ldots \mathsf{R}_{M_n}(k_n). 
\end{align}
Let us now make a substitution 
\begin{align} 
m_1 &:= k_1 + \frac{B_{M_1}}{2 A_{M_1}}, \ldots, m_n := k_n + \frac{B_{M_n}}{2 A_{M_n}}.
\end{align}
Let $\mathsf{Q}_{M_i}(t)$ be as defined by \eqref{QM_defn}, i.e. 
\begin{align} 
\mathsf{Q}_{M_i}(t)&:= \mathsf{R}_{M_i}\Big(t - \frac{B_{M_i}}{2 A_{M_i}} \Big).
\end{align}
With this, the counting function becomes 
\begin{align}
N_M(\lambda) & = \sum_{\substack{m_1-\frac{B_{M_1}}{2A_{M_1}}, \ldots, m_n-\frac{B_{M_n}}{2A_{M_n}} 
\in \mathbb{Z}^{\geq 0}\\ A_{M_{1}} m_1^2 + \ldots + A_{M_{n}} m_n^2 \leq \lambda^2 + \kappa}}
\mathsf{Q}_{M_1}(m_1) \ldots \mathsf{Q}_{M_n}(m_n), 
\end{align}
where $\kappa$ is given by 
\begin{align*}
\kappa & := \sum_{i=1}^{n}\Big(\frac{B_{M_i}^2}{4 A_{M_{i}}} - C_{M_i}\Big). 
\end{align*}
Since each of the $M_i$ is a 
$\mathsf{W}$-manifold, we can use theorem \ref{Auxiliary_result} 
and conclude that 
\begin{align}
N_M(\lambda) & = C\lambda^d + O(\lambda^{d-1-\frac{n-1}{n+1}}).
\end{align}
The constant $C$ can now be obtained from Weyl's law for compact manifolds without boundary.  

\section{CROSSes are $\mathsf{W}$-manifolds: proof} 
\label{CROSS_is_W}
We will now prove theorem \ref{CROSS_is_W_ag}, i.e. we will 
verify that all the CROSSes are $\mathsf{W}$-manifolds.\\ 
\hf \hf First of all, we recall that the eigenvalues and multiplicities of the 
Laplacian can be packaged together in the zeta function. 
If we are in the setting of section \ref{setup_not}, we get that 
\begin{align}
\zeta_{M}(s) & = \sum_{k=0}^{\infty} \mathsf{R}_M(k) \exp(-\alpha_M(k)s). \label{zeta_general}
\end{align}  
We will now use the results of the paper \cite{Cahn}, where the authors explicitly compute 
$\zeta_{M}(s)$ for all CROSSes. Knowing $\zeta_{M}(s)$ immediately tells us what is 
$\alpha_M(k)$ and $\mathsf{R}_M(k)$. Let us now explicitly write down the zeta function for 
each of these spaces. Although the case of the sphere has been worked out in \cite{Weyl_Wyman_Iosevich}, we will  
include it here for completeness.  

\subsection{Odd dimensional spheres}
Let us start with the odd dimensional sphere $M=S^{2d-1}$, with $d \geq 3$. Following 
\cite[page 7, theorem 2.4]{Cahn}, we conclude that 
\begin{align}
\zeta_{M}(s) & = \sum_{k=0}^{\infty} \Big(\frac{k+d-1}{d-1}\prod_{j=1}^{2d-3} \frac{k+j}{j}\Big) 
\exp\Big(-\Big(\frac{k^2+2k(d-1)}{4(d-1)}\Big)s\Big). \label{seta_S_d}
\end{align}
Comparing \eqref{seta_S_d} with \eqref{zeta_general}, we conclude that 
\begin{align}
\alpha_{M}(k)& = \Big(\frac{k^2+2k(d-1)}{4(d-1)}\Big) \qquad \textnormal{and} \qquad 
\mathsf{R}_{M}(k) = 
\Big(\frac{k+d-1}{d-1}\prod_{j=1}^{2d-3} \frac{k+j}{j}\Big).
\end{align}
Note that $\mathsf{R}_{M}(k)$ is a polynomial of degree $2d-2$. Now, 
\[\alpha_{M}(k)=\frac{1}{4(d-1)}k^{2}+\frac{1}{2}k=A_{M}k^{2}+B_{M}k+C_{M}, \]
where $A_{M}=\frac{1}{4(d-1)}, B_{M}=\frac{1}{2}, C_{M}=0 $ and therefore, 
\[\frac{B_{M}}{2 A_{M}}=d-1. \]
Next, we apply definition \eqref{QM_defn} and conclude that 
\begin{align}
\mathsf{Q}_{M}(t) & = \mathsf{R}_{M}(t-(d-1))=\frac{t}{d-1} \prod_{j=1}^{2d-3} \frac{t-(d-1)+j}{j}=\frac{t}{d-1} \prod_{j=1}^{2d-3} \frac{t+a_{j}}{j} ,
\end{align}
where $a_{j}:=j-(d-1),\ j=1,\dots,2d-3 $. \\ 
\hf \hf We claim that the coefficient of $t^{2d-3} $ in $\mathsf{Q}_{M}(t)$ is zero. To see this, we note that it 
is sufficient to prove that the numerator of the coefficient of $t^{2d-4} $ in $ \prod_{j=1}^{2d-3} \frac{t+a_{j}}{j} $, that is, $\sum_{j=1}^{2d-3} a_{j} $ is zero. Now,
\begin{equation}
\begin{aligned}
\sum_{j=1}^{2d-3} a_{j}&=\sum_{j=1}^{2d-3} \left(j-d+1 \right)=\frac{(2d-3)(2d-2)}{2}-(d-1)(2d-3)=0,
\end{aligned}
\notag
\end{equation}
and hence, we are done. \\ 
\hf \hf To finish the case of odd spheres, we also need to separately consider $M=S^3$. 
Following 
\cite[page 6, equation 2.3 (a)]{Cahn}, we conclude that 
\begin{align}
\zeta_{M}(s) & = \sum_{k=0}^{\infty} (k+1)^2 \exp\Big(-\Big(\frac{k^2+2k}{4}\Big)s\Big). \label{seta_S_3}
\end{align}
Comparing \eqref{seta_S_3} with \eqref{zeta_general}, we conclude that 
\begin{align}
\alpha_{M}(k)& = \Big(\frac{k^2+2k}{4}\Big) \qquad \textnormal{and} \qquad 
\mathsf{R}_{M}(k) = (k+1)^2.  
\end{align}
Note that $\mathsf{R}_{M}(k)$ is a polynomial of degree $2$. Now, 
\[\alpha_{M}(k)=\frac{1}{4}k^{2}+\frac{1}{2}k=A_{M}k^{2}+B_{M}k+C_{M}, \]
where $A_{M}=\frac{1}{4}, B_{M}=\frac{1}{2}, C_{M}=0 $ and therefore, 
\[\frac{B_{M}}{2 A_{M}}=1. \] 
Next, we apply definition \eqref{QM_defn} and conclude that the coefficient of $t$ in the polynomial
\begin{align}
\mathsf{Q}_{M}(t) & = \mathsf{R}_{M}(t-1)
\end{align}
is zero.
This completes all the odd spheres.

\subsection{Even dimensional spheres}
 Next, we consider even dimensional spheres $M=S^{2d}$, with $d\geq 1$.  
Following 
\cite[page 8, theorem 3.2]{Cahn}, we conclude that 
\begin{align}
\zeta_{M}(s) & = \sum_{k=0}^{\infty} \Big(\frac{2k+2d-1}{2d-1}\prod_{j=1}^{2d-2} \frac{k+j}{j}\Big) 
\exp\Big(-\Big(\frac{k^2+k(2d-1)}{4d-2}\Big)s\Big). \label{seta_S_d_even}
\end{align}
Comparing \eqref{seta_S_d_even} with \eqref{zeta_general}, we conclude that 
\begin{align}
\alpha_{M}(k)& = \Big(\frac{k^2+k(2d-1)}{4d-2}\Big) \qquad \textnormal{and} \qquad 
\mathsf{R}_{M}(k) = 
\Big(\frac{2k+2d-1}{2d-1}\prod_{j=1}^{2d-2} \frac{k+j}{j}\Big).
\end{align}
Note that $\mathsf{R}_{M}(k)$ is a polynomial of degree $2d-1$.  Now, 
\[\alpha_{M}(k)=\frac{1}{4d-2}k^{2}+\frac{2d-1}{4d-2}k=A_{M}k^{2}+B_{M}k+C_{M}, \]
where $A_{M}=\frac{1}{4d-2}, B_{M}=\frac{1}{2}, C_{M}=0 $ and therefore, 
\[\frac{B_{M}}{2 A_{M}}=\frac{2d-1}{2}. \] 
Next, we apply definition \eqref{QM_defn} and conclude that 
\begin{align}
\mathsf{Q}_{M}(t) & = \mathsf{R}_{M}\Big(t-\frac{2d-1}{2}\Big)= \frac{2}{(2d-1)}t \cdot \prod_{j=1}^{2d-2} \frac{t+a_{j}}{j},
\end{align}
where $a_{j}:=-\frac{2d-1}{2}+j,\ j=1,\dots,2d-2 $. \\
\hf \hf We claim that the coefficient of $t^{2d-2} $ in $\mathsf{Q}_{M}(t)$ is zero. To see this, we note that it is sufficient to prove that the numerator of the coefficient of $t^{2d-3} $ in $ \prod_{j=1}^{2d-2} \frac{t+a_{j}}{j} $, that is, $\sum_{j=1}^{2d-2} a_{j} $ is zero. Now,
\begin{equation}
\begin{aligned}
\sum_{j=1}^{2d-2} a_{j}&=\sum_{j=1}^{2d-2} \left(j-\frac{2d-1}{2} \right)=\frac{(2d-2)(2d-1)}{2}-\frac{(2d-1)(2d-2)}{2}=0,
\end{aligned}
\notag
\end{equation}
and hence, we are done. This completes all the spheres. 
%
%
\subsection{Odd dimensional real projective spaces}
Let us now move on to the real projective spaces. 
We will start with the odd dimensional real projective spaces 
$M=\mathbb{RP}^{2d-1}$, with $d \geq 3$. 
Following \cite[page 8, corollary 2.6]{Cahn}, we conclude that 
\begin{align}
\zeta_{M}(s) & = \sum_{k=0}^{\infty} \Big(\frac{2k+d-1}{d-1}\prod_{j=1}^{2d-3} \frac{2k+j}{j}\Big) 
\exp\Big(-\Big(\frac{k^2+k(d-1)}{d-1}\Big)s\Big). \label{seta_RPd}
\end{align}
Comparing \eqref{seta_RPd} with \eqref{zeta_general}, we conclude that 
\begin{align}
\alpha_{M}(k)& = \left(\frac{k^2+k(d-1)}{d-1}\right) \qquad \textnormal{and} \qquad 
\mathsf{R}_{M}(k) = 
\left(\frac{2k+d-1}{d-1}\prod_{j=1}^{2d-3} \frac{2k+j}{j}\right).
\end{align}
Note that $\mathsf{R}_{M}(k)$ is a polynomial of degree $2d-2$.  Now, 
\[\alpha_{M}(k)=\frac{1}{d-1}k^{2}+k=A_{M}k^{2}+B_{M}k+C_{M}, \]
where $A_{M}=\frac{1}{d-1}, B_{M}=1, C_{M}=0 $ and therefore, 
\[\frac{B_{M}}{2 A_{M}}=\frac{d-1}{2}. \] 
Next, we apply definition \eqref{QM_defn} and conclude that 
\begin{align}
\mathsf{Q}_{M}(t) & = \mathsf{R}_{M}\Big(t-\frac{d-1}{2}\Big)=\frac{2t}{d-1} \prod_{j=1}^{2d-3}  \frac{2t+a_{j}}{j},
\end{align}
where $a_{j}:=j-(d-1),\ j=1,\dots,2d-3$. \\
\hf \hf We claim that the coefficient of $t^{2d-3} $ in $\mathsf{Q}_{M}(t)$ is zero. To see this, 
we note that it is sufficient to prove that the numerator of the coefficient of 
$t^{2d-4} $ in $ \prod_{j=1}^{2d-3} \frac{2t+a_{j}}{j}=\prod_{j=1}^{2d-3} \frac{t+\frac{a_{j}}{2}}{\frac{j}{2}} $, that is, $\sum_{j=1}^{2d-3} \frac{a_{j}}{2} $ is zero. Now,
\begin{equation}
\begin{aligned}
\sum_{j=1}^{2d-3} \frac{a_{j}}{2}&=\frac{1}{2}\sum_{j=1}^{2d-3} \left(j-d+1 \right)=\frac{1}{2} \left(\frac{(2d-3)(2d-2)}{2}-(d-1)(2d-3) \right)=0,
\end{aligned}
\notag
\end{equation}
and hence, we are done. \\ 
\hf \hf To finish the case of the odd real projective spaces, 
we also need to separately consider $M=\mathbb{RP}^3$. 
Following 
\cite[page 6, equation 2.3 (b)]{Cahn}, we conclude that 
\begin{align}
\zeta_{M}(s) & = \sum_{k=0}^{\infty} (2k+1)^2 \exp\Big(-(k^2+k)s\Big). \label{zeta_RP3}
\end{align}
Comparing \eqref{zeta_RP3} with \eqref{zeta_general}, we conclude that 
\begin{align}
\alpha_{M}(k)& = k^2 + k\qquad \textnormal{and} \qquad 
\mathsf{R}_{M}(k) = (2k+1)^2.  
\end{align}
Note that $\mathsf{R}_{M}(k)$ is a polynomial of degree $2$. Now, 
\[\alpha_{M}(k)=A_{M}k^{2}+B_{M}k+C_{M}, \]
where $A_{M}=1, B_{M}=1, C_{M}=0 $ and therefore, 
\[\frac{B_{M}}{2 A_{M}}=\frac{1}{2}. \] 
Next, we apply definition \eqref{QM_defn} and observe that the coefficient of $t$ in the polynomial
\begin{align}
\mathsf{Q}_{M}(t) & = \mathsf{R}_{M}\Big(t-\frac{1}{2}\Big)=4t^{2}
\end{align}
is zero. 
This completes all the odd real projective spaces.  

\subsection{Even dimensional real projective spaces}
Next, we consider even dimensional real projective space $M=\mathbb{RP}^{2d}$, with $d\geq 1$.  
Following 
\cite[pages 8-9, corollary 3.4]{Cahn}, we conclude that 
\begin{align}
\zeta_{M}(s) & = \sum_{k=0}^{\infty} \Big(\frac{4k+2d-1}{2d-1}\prod_{j=1}^{2d-2} \frac{2k+j}{j}\Big) 
\exp\Big(-\Big(\frac{2k^2+k(2d-1)}{2d-1}\Big)s\Big). \label{seta_RPd_even}
\end{align}
Comparing \eqref{seta_RPd_even} with \eqref{zeta_general}, we conclude that 
\begin{align}
\alpha_{M}(k)& = \Big(\frac{2k^2+k(2d-1)}{2d-1}\Big) \qquad \textnormal{and} \qquad 
\mathsf{R}_{M}(k) = 
\Big(\frac{4k+2d-1}{2d-1}\prod_{j=1}^{2d-2} \frac{2k+j}{j}\Big).
\end{align}
Note that $\mathsf{R}_{M}(k)$ 
is a polynomial of degree $2d-1$.  Now, 
\[\alpha_{M}(k)=\frac{2}{2d-1}k^{2}+k=A_{M}k^{2}+B_{M}k+C_{M}, \]
where $A_{M}=\frac{2}{2d-1}, B_{M}=1, C_{M}=0 $ and therefore, 
\[\frac{B_{M}}{2 A_{M}}=\frac{2d-1}{4}. \] 
Next, we apply definition \eqref{QM_defn} and conclude that 
\begin{align}
\mathsf{Q}_{M}(t) & = \mathsf{R}_{M}\Big(t-\frac{2d-1}{4}\Big)=\frac{4t}{2d-1}  \prod_{j=1}^{2d-2} \frac{2t+a_{j}}{j} ,
\end{align}
where $a_{j}:=j-\frac{2d-1}{2},\ j=1,\dots,2d-2$. We claim that the coefficient of $t^{2d-2} $ in $\mathsf{Q}_{M}(t)$ is zero. 
To see this, we note that it is sufficient to prove that the numerator of the coefficient of $t^{2d-3} $ in $ \prod_{j=1}^{2d-2} \frac{2t+a_{j}}{j} =\prod_{j=1}^{2d-2} \frac{t+\frac{a_{j}}{2}}{\frac{j}{2}}$, that is, $\sum_{j=1}^{2d-2} \frac{a_{j}}{2} $ is zero. Now,
\begin{equation}
\begin{aligned}
\sum_{j=1}^{2d-2} \frac{a_{j}}{2}&=\frac{1}{2}\sum_{j=1}^{2d-2} \left(j-\frac{2d-1}{2} \right)=\frac{1}{2} \left(\frac{(2d-2)(2d-1)}{2}-\frac{(2d-1)(2d-2)}{2} \right)=0,
\end{aligned}
\notag
\end{equation}
and hence, we are done.

\subsection{Complex projective spaces}
Let us now move on to the complex projective spaces $M=\mathbb{CP}^d$, 
with $d \geq 2$. Notice that the real dimension of this manifold is $2d$. 
Following \cite[page 9, theorem 4.2]{Cahn}, we conclude that 
\begin{align}
\zeta_{M}(s) & = \sum_{k=0}^{\infty} \Big(\frac{2k+d}{d}\prod_{j=1}^{d-1} \Big(\frac{k+j}{j}\Big)^2\Big) 
\exp\Big(-\Big(\frac{k^2+kd}{d+1}\Big)s\Big). \label{zeta_CPd}
\end{align}
Comparing \eqref{zeta_CPd} with \eqref{zeta_general}, we conclude that 
\begin{align}
\alpha_{M}(k)& = \left(\frac{k^2+kd}{d+1}\right) \qquad \textnormal{and} \qquad 
\mathsf{R}_{M}(k) = 
\left(\frac{2k+d}{d}\prod_{j=1}^{d-1} \Big(\frac{k+j}{j}\Big)^2\right).
\end{align}
Note that $\mathsf{R}_{M}(k)$ is a polynomial of degree $2d-1$. Now, 
\[\alpha_{M}(k)=\frac{1}{d+1}k^{2}+\frac{d}{d+1}k=A_{M}k^{2}+B_{M}k+C_{M}, \]
where $A_{M}=\frac{1}{d+1}, B_{M}=\frac{d}{d+1}, C_{M}=0 $ and therefore, 
\[\frac{B_{M}}{2 A_{M}}=\frac{d}{2}. \] 
Next, we apply definition \eqref{QM_defn} and conclude that 
\begin{align}
\mathsf{Q}_{M}(t) & = \mathsf{R}_{M}\left(t-\frac{d}{2}\right)= \frac{2t}{d} \prod_{j=1}^{d-1} \left(\frac{t+a_{j}}{j} \right)^{2},
\end{align}
where $a_{j}:= j-\frac{d}{2},\ j=1,\dots, d-1 $.  We claim that the coefficient of $t^{2d-2} $ in $\mathsf{Q}_{M}(t)$ is zero. To see this, we note that it is sufficient to prove that the numerator of the coefficient of $t^{2d-3} $ in $\prod_{j=1}^{d-1} \left(\frac{t+a_{j}}{j} \right)^{2} $, that is, $\sum_{j=1}^{d-1} 2a_{j} $ is zero. Now,
\begin{equation}
\begin{aligned}
\sum_{j=1}^{d-1} 2a_{j}&=2\sum_{j=1}^{d-1} \left(j-\frac{d}{2} \right)=2 \left(\frac{d(d-1)}{2}-\frac{d(d-1)}{2} \right)=0,
\end{aligned}
\notag
\end{equation}
and hence, we are done.
\subsection{Quaternionic projective spaces}
Let us now move on to the quaternionic projective spaces $M=\mathbb{HP}^{d-1}$, 
with $d \geq 2$. 
We are referring to the (quaternionic) dimension as $d-1$ as opposed to $d$ in order to make it easy to follow the 
notation of \cite{Cahn}. 
Notice that the real dimension of 
$\mathbb{HP}^{d-1}$
is $4d-4$. 
Following \cite[page 10, theorem 5.2]{Cahn}, we conclude that 
\begin{align}
\zeta_{M}(s) & = \sum_{k=0}^{\infty} \Big(\frac{2k+2d-1}{2d-1}\prod_{j=2}^{2d-2} \frac{k+j}{j} 
\prod_{l=1}^{2d-3} \frac{k+l}{l}  \Big)
\exp\Big(-\Big(\frac{k^2+2kd-k}{2(d+1)}\Big)s\Big). \label{zeta_HPd}
\end{align}
Comparing \eqref{zeta_HPd} with \eqref{zeta_general}, we conclude that 
\begin{align}
\alpha_{M}(k)& = \Big(\frac{k^2+2kd-k}{2(d+1)}\Big) \qquad \textnormal{and} \nonumber \\ 
\mathsf{R}_{M}(k) & = 
\Big(\frac{2k+2d-1}{2d-1}\prod_{j=2}^{2d-2} \frac{k+j}{j} 
\prod_{l=1}^{2d-3} \frac{k+l}{l}  \Big). \nonumber 
\end{align}
Notice that $\mathsf{R}_{M}(k)$ is a polynomial of degree $4d-5$. Now, 
\[\alpha_{M}(k)=\frac{1}{2(d+1)}k^{2}+\frac{2d-1}{2(d+1)}k=A_{M}k^{2}+B_{M}k+C_{M}, \]
where $A_{M}=\frac{1}{2(d+1)}, B_{M}=\frac{2d-1}{2(d+1)}, C_{M}=0 $ and therefore, 
\[\frac{B_{M}}{2 A_{M}}=\frac{2d-1}{2}=d-\frac{1}{2}. \] 
Next, we apply definition \eqref{QM_defn} and conclude that 
\begin{align}
\mathsf{Q}_{M}(t) & = \mathsf{R}_{M}\left(t-\left(d-\frac{1}{2} \right)\right)=\frac{2t}{2d-1} \prod_{j=2}^{2d-2} \frac{t+a_{j}}{j} \prod_{l=1}^{2d-3} \frac{t+b_{l}}{l},
\end{align}
where 
\begin{equation}
\begin{aligned}
a_{j}&= j-\left(d-\frac{1}{2} \right), \ j=2,\dots,2d-2, \\
b_{l}&= l-\left(d-\frac{1}{2} \right), \ l=1,\dots,2d-3.
\end{aligned}
\notag
\end{equation}
We claim that the coefficient of $t^{4d-6} $ in $\mathsf{Q}_{M}(t)$ is zero. To see this, we note that it is sufficient to prove that the numerator of the coefficient of $t^{4d-7} $ in $  \prod_{j=2}^{2d-2} \frac{t+a_{j}}{j} \prod_{l=1}^{2d-3} \frac{t+b_{l}}{l} $, that is, $\sum_{j=2}^{2d-2} a_{j}+\sum_{l=1}^{2d-3} b_{l} $ is zero. Now,
\begin{equation}
\begin{aligned}
\sum_{j=2}^{2d-2} a_{j}+\sum_{l=1}^{2d-3} b_{l}&=\sum_{j=2}^{2d-2} j+\sum_{l=1}^{2d-3} l-(4d-6)\left(d-\frac{1}{2} \right)\\
&=\left(\frac{(2d-2)(2d-1)}{2}-1\right)+\frac{(2d-3)(2d-2)}{2}-\frac{(4d-6)(2d-1)}{2}=0,
\end{aligned}
\notag
\end{equation}
and hence, the claim is true.
\subsection{Octonionic projective plane}
Finally, we consider $M=\mathbb{OP}^2$, the octonionic projective plane.  
Note that the real dimension of this space is $16$. 
Following 
\cite[page 11, theorem 6.2]{Cahn}, we conclude that 
\begin{align}
\zeta_{M}(s) & = \sum_{k=0}^{\infty} \frac{2k+11}{11} 
\left(\prod_{j=1}^{3} \frac{k+j}{j}\right) \left(\prod_{l=4}^{7} \left(\frac{k+l}{l}\right)^2\right) 
\left(\prod_{m=8}^{10} \frac{k+m}{m}\right)  \exp\left(-\frac{(k^2+11k)}{18}s\right). \label{zeta_OP2}
\end{align}
Comparing \eqref{zeta_OP2} with \eqref{zeta_general}, we conclude that 
\begin{align}
\alpha_{M}(k)& = \frac{(k^2+11k)}{18} \qquad \textnormal{and} \nonumber \\ 
\mathsf{R}_{M}(k) & = 
\frac{2k+11}{11} 
\left(\prod_{j=1}^{3} \frac{k+j}{j}\right) \left(\prod_{l=4}^{7} \left(\frac{k+l}{l}\right)^2\right) 
\left(\prod_{m=8}^{10} \frac{k+m}{m}\right).  \nonumber 
\end{align}
Note that $\mathsf{R}_{M}(k)$ is a polynomial of degree $15$.  Now, 
\[\alpha_{M}(k)=\frac{1}{18}k^{2}+\frac{11}{18}k=A_{M}k^{2}+B_{M}k+C_{M}, \]
where $A_{M}=\frac{1}{18}, B_{M}=\frac{11}{18}, C_{M}=0 $ and therefore, 
\[\frac{B_{M}}{2 A_{M}}=\frac{11}{2}. \] 
Next, we apply definition \eqref{QM_defn} and conclude that 
\begin{align}
\mathsf{Q}_{M}(t) & = \mathsf{R}_{M}\Big(t-\frac{11}{2}\Big)= \frac{2t}{11} \prod_{j=1}^{3} \frac{t+a_{j}}{j} \prod_{l=4}^{7} \left( \frac{t+a_{l}}{l} \right)^{2} \prod_{m=8}^{10} \frac{t+a_{m}}{m},
\end{align}
where
\begin{equation}
\begin{aligned}
a_{j}&= j-\frac{11}{2}, \ j=1,\dots,10.
\end{aligned}
\notag
\end{equation}
We claim that the coefficient of $t^{14} $ in $\mathsf{Q}_{M}(t)$ is zero. To see this, we note that it is sufficient to prove that the numerator of the coefficient of $t^{13} $ in $ \prod_{j=1}^{3} \frac{t+a_{j}}{j} \prod_{l=4}^{7} \left( \frac{t+a_{l}}{l} \right)^{2} \prod_{m=8}^{10} \frac{t+a_{m}}{m} $, that is, $\sum_{j=1}^{3} a_{j}+\sum_{l=4}^{7} 2a_{l} +\sum_{m=8}^{10} a_{m} $ is zero. Now,
\begin{equation}
\begin{aligned}
\sum_{j=1}^{3} a_{j}+\sum_{l=4}^{7} 2a_{l} +\sum_{m=8}^{10} a_{m}&=\sum_{j=1}^{3} j+\sum_{m=8}^{10} m+2 \sum_{l=4}^{7} l+ 3 \left( -\frac{11}{2}\right)+ 3  \left( -\frac{11}{2}\right)+8  \left( -\frac{11}{2}\right)=0,
\end{aligned}
\notag
\end{equation}
and hence, the claim is proved.\\
This completes the proof of theorem \ref{CROSS_is_W_ag}. \qed 
 

\section{Error is sharp for a single $\mathsf{W}$-manifold: proof} 
\label{Error_is_sharp}
In this section, we will prove theorem \ref{main_thm_W2_sharp}.
It was shown by H{\"o}rmander in \cite{LH} that the error term in Weyl's law can not 
be improved for the round sphere.  
We will generalize his proof to $\mathsf{W}$-manifolds. 
Let us first prove a few auxiliary results.  

\begin{prp}
\label{psi_k}
Let us define the following function
$\psi_k: \mathbb{N}\longrightarrow \mathbb{R}$, given by   
\begin{align} 
\psi_k(n)&:= \sum_{i=1}^{n} i^{k-1} \nonumber \\ 
          & = 1^{k-1} + 2^{k-1} + \ldots + n^{k-1},   
\end{align} 
where $k \in \mathbb{Z}^{\geq 2}$. Then, 
\begin{align} 
\psi_k(n)& = \frac{n^k}{k} + \frac{1}{2} n^{k-1} + \mathcal{R}_{k-2}(n), \label{nk_error_sharp}
\end{align}
where $\mathcal{R}_{k-2}(n)$ is a polynomial in $n$ of degree at most $k-2$. 
\end{prp}
\begin{proof}
Let $n \in \mathbb{N} $ and $k\in \mathbb{Z}^{\geq 2} $. 
Note that
\begin{equation}
\label{tel_s}
\begin{aligned}
&(n+1)^{k}-n^{k}= kn^{k-1}+\frac{k(k-1)}{2} n^{k-2}+ \cdots + 1\\
&n^{k}-(n-1)^{k}= k(n-1)^{k-1}+\frac{k(k-1)}{2} (n-1)^{k-2}+ \cdots + 1\\
&\cdots \\
&\cdots \\
&2^{k}-1^{k}= k \cdot 1^{k-1}+ \frac{k(k-1)}{2} \cdot 1^{k-2} + \cdots+ 1.
\end{aligned}
\end{equation}
Adding up all the terms on 
the left hand side of 
equation 
\eqref{tel_s}, we get a telescopic sum which equals 
$(n+1)^{k}-1^{k}$. 
Equating this with the summation of the right hand 
side of equation \eqref{tel_s}, 
we conclude that     
\begin{equation}
(n+1)^{k}-1^{k} = k\psi_{k}(n)+\frac{k(k-1)}{2} \psi_{k-1}(n)+\cdots+\psi_1(n). 
\label{nk_sum}
\end{equation}
In equation \eqref{nk_sum}, we note that $\psi_{1}(n)=n $.
Using the fact that $\psi_2(n) = \frac{n(n+1)}{2}$ 
and equation \eqref{nk_sum}, equation \eqref{nk_error_sharp} 
follows by using induction on $k$. 
This proves the desired result.  
\end{proof}
Next, we prove the following result. 
\begin{lmm}
\label{sharp_main}
Let $a \in \mathbb{R}_{>0}$, $b\in \mathbb{R}$ and $y \in \mathbb{R}_{\geq 0}$. For $\lambda$ 
sufficiently large and positive, define 
\begin{align*}
M(\lambda)& = \frac{[\sqrt{a\lambda^2+b}-y]^d}{d} + (y+\frac{1}{2})[\sqrt{a\lambda^2+b}-y]^{d-1}.
\end{align*}
Then
\begin{align}
\lim_{\lambda\longrightarrow +\infty} \frac{M(\lambda)}{\lambda^d} & = \frac{a^{\frac{d}{2}}}{d}. 
\label{M_leading}
\end{align}
Furthermore, the error term $\mathcal{E}(\lambda)$ can not be $O(\lambda^{d-1-\varepsilon})$ 
for any positive $\varepsilon$, where 
\begin{align*}
\mathcal{E}(\lambda) &:= M(\lambda)-\frac{a^{\frac{d}{2}}}{d} \lambda^d.
\end{align*}
\end{lmm} 
\noindent \textbf{Proof:} Equation \eqref{M_leading} is immediate using the fact that 
\begin{align} 
A+B-1&\leq [A+B]\leq A+B+1. \label{A_B_box}
\end{align}
An application of \eqref{A_B_box}, binomial theorem 
and the squeeze rule immediately gives us \eqref{M_leading}.  
We now prove the second assertion. We will prove by contradiction. Suppose 
$\mathcal{E}(\lambda)$ was in $O(\lambda^{d-1-\varepsilon})$ for some $\varepsilon>0$. 
Let us define the following discrete subset of $\mathbb{R}$ given by 
\begin{align*}
\Lambda&:= \Big\{ \lambda\in \mathbb{R}_{>0}: a \lambda^2+b>0, \sqrt{a\lambda^2+b}-y \in \mathbb{Z}_{\geq 0}
\Big\}.
\end{align*}
Define $\gamma:\Lambda\longrightarrow \mathbb{Z}_{>0}$, given by 
\begin{align} 
\gamma(\lambda)&:= \sqrt{a\lambda^2+b}-y. \label{gamma_defn_in_terms_of_lambda} 
\end{align}
To lighten the notation, 
we will abbreviate $\gamma(\lambda)$ by $\gamma$; however, it should always 
be understood that $\gamma$ depends on $\lambda$ via equation \eqref{gamma_defn_in_terms_of_lambda}. 
Hence, if the error term was in $O(\lambda^{d-1-\varepsilon})$, then it would imply that there is some constant 
$C>0$ such that for all $\lambda \in \Lambda$, sufficiently large, we have 
\begin{align}
\Big| \Big(y+\frac{1}{2}\Big) \gamma^{d-1} 
+\Big(\frac{\gamma^d}{d}-\frac{a^{\frac{d}{2}}}{d} \lambda^d\Big) \Big| & \leq C \lambda^{d-1-\varepsilon}. \label{lm_gmm_contrd}
\end{align} 
Using the fact that $|A|-|B|\leq |A+B|$
(note that $|A+B-B| \leq |A+B| + |B|$ which implies what we have written), 
we conclude that 
\begin{align}
\Big|\Big(y+\frac{1}{2}\Big) \gamma^{d-1}\Big| 
-\Big| \frac{\gamma^d}{d}-\frac{a^{\frac{d}{2}}}{d} \lambda^d \Big| 
& \leq C \lambda^{d-1-\varepsilon}. \label{lm_gmm_contrd_ag}
\end{align} 
Here we set 
\[A:= \Big(y+\frac{1}{2}\Big) \gamma^{d-1} \qquad \textnormal{and} 
\qquad B:= \Big(\frac{\gamma^d}{d}-\frac{a^{\frac{d}{2}}}{d} \lambda^d\Big). \]
Now we divide both sides of equation \eqref{lm_gmm_contrd_ag} 
and take the limit as $\lambda$ tends to infinity (and $\lambda \in \Lambda$). 
The right hand side will obviously be zero. We first note that 
\begin{align*} 
\lim_{\substack{\lambda \longrightarrow \infty\\ \lambda \in \Lambda}} 
\frac{\gamma^d}{d} - \frac{a^{\frac{d}{2}}}{d} \lambda^d & = 0. 
\end{align*}
This follows by taking a binomial expansion of $\gamma$ and noticing that after subtracting the 
leading term (namely $\frac{a^{\frac{d}{2}}}{d} \lambda^{d}$), we are left with a terms of order $\lambda^{d-2}$ or less. 
Hence, what remains is the first term on the left hand side of equation 
\eqref{lm_gmm_contrd_ag}. 
We note that 
\begin{align*} 
\lim_{\substack{\lambda \longrightarrow \infty\\ \lambda \in \Lambda}} 
\Big(y+\frac{1}{2}\Big) \frac{\gamma^{d-1}}{\lambda^{d-1}} & = \Big(y+ \frac{1}{2}\Big).
\end{align*} 
Since $y$ is non negative, the above limit is positive. This gives us a contradiction. \qed

\begin{thm}
\label{Auxiliary_result_2}
Let $d \in \mathbb{Z}^{\geq 2}$  
and $\mathsf{Q}_{d-1}(t)$ be a degree $d-1$ polynomial in $t$ that does not have the $t^{d-2}$ term. In other words, 
\begin{align}
\mathsf{Q}_{d-1}(t) & = C t^{d-1} + \varphi(t), \qquad \textnormal{if} \qquad d>2 \qquad \textnormal{and} \nonumber \\ 
\mathsf{Q}_{d-1}(t) & = C t, \qquad \textnormal{if} \qquad d=2, \nonumber 
\end{align}
where $\varphi(t)$ is a polynomial of degree at most $d-3$ and 
$C$ is some non-zero constant. 
Let $y$ be a fixed non negative real number and 
$A$ a fixed positive real numbers. 
Define 
$N:\mathbb{R}\longrightarrow \mathbb{R}$ 
as  
\begin{align}
N(\lambda) & := \sum_{\substack{m-y \in \mathbb{Z}^{\geq 0}\\ A m^2 \leq a \lambda^2 + b}}
\mathsf{Q}_{d-1}(m). 
\end{align} 
where $a \in \mathbb{R}^{>0}$ and $b \in \mathbb{R}$. 
Then the error term in the asymptotic expression of $N(\lambda)$ 
(after subtracting off the leading term) can not be improved by $O(\lambda^{d-1-\varepsilon})$ 
for any positive $\varepsilon$. 
\end{thm} 

\begin{proof} Without loss of generality, let us set $C$ to be equal to one.  
We sum up the leading term of $\mathsf{Q}_{d-1}(t)$ 
(namely the degree $t^{d-1}$ terms) using Proposition \ref{psi_k}. 
The first two terms precisely correspond to $M(\lambda)$ as given by 
Lemma \ref{sharp_main} (replace $a$ by $\frac{a}{A}$ and replace $b$ by $\frac{b}{A}$). 
The remaining terms are all of order $\lambda^{d-2}$ or less 
(this is where we use the fact that $\mathsf{Q}_{d-1}(t)$ does not have a $t^{d-2}$ term). 
Hence the remaining terms can not compensate 
for the error term from $M(\lambda)$. We have already seen that the error term for $M(\lambda)$ 
can not be improved. This shows that the error term for $N(\lambda)$ can not be improved.  
\end{proof}

Since the Eigenvalue counting function $N(\lambda)$ 
for all $\mathsf{W}$-manifolds is of the form 
considered in Theorem \ref{Auxiliary_result_2}, we have proved Theorem \ref{MR2} 

	
\section{Counting Lattice Points in a Ball}
\label{ball_std}
To keep the paper self contained, 
we give a proof of  
equation \eqref{M_lm} (even though it is a standard fact). 
We will prove the following general result.  
\begin{prp}
\label{lattice_count}
Let $p$ be any norm on $\mathbb{R}^d$ where $d \geq 1$. Given $\lambda > 0$,  let $B_p(\lambda)$ denote the closed $p$-ball of radius $\lambda$ around the origin, i.e.  \(B_p(\lambda) := \{x \in \mathbb{R}^d : p(x) \leq \lambda\}\). 
Let \(N^p(\lambda) := \#(\mathbb{Z}^d \cap B_p(\lambda))\) denote 
the number of integer lattice points in $B_p(\lambda)$.  
Then, $N^p(\lambda) \in O(\lambda^d)$ as $\lambda \to \infty$.
\end{prp}

The supremum norm on $\mathbb{R}^d$ is particularly suited for counting integer lattice points, so this is the case we first deal with.  
\begin{lmm}
	\label{sup_norm}
	Let $p_{\infty}$ denote the supremum norm on $\mathbb{R}^d$, i.e. \(p_{\infty}(x) := \max_{1 \leq i \leq n}|x_i|\). Then $\forall \lambda > 0$, $N^{p_{\infty}}(\lambda) \leq (2\lambda + 3)^d $.
\end{lmm}
We will simplify notation slightly:  $B_{\infty}(\lambda)$ will denote $B_{p_{\infty}}(\lambda)$ and $N^{\infty}(\lambda)$ will denote $N^{p_{\infty}}(\lambda)$.
\begin{proof}
	Observe that $B_{\infty}(\lambda) = [-\lambda,\lambda]^d$ is a $d$-dimensional cube centered at the origin. \\We will first assume $\lambda \in \mathbb{N}$. It is then clear that  \[x \in \mathbb{Z}^d \cap B_{\infty}(\lambda) \iff x_i \in \mathbb{Z} \cap [-\lambda,\lambda] \text{ for each } 1 \leq i \leq d\]. \\
	As $\lambda \in \mathbb{N}$, this means each $x_i$ can take one of precisely $2\lambda + 1$ values, i.e \[x_i \in \{-\lambda, -\lambda +1, \ldots, 0, \ldots ,\lambda-1, \lambda\}\]
	and we conclude \begin{equation}
		N^{\infty}(\lambda) = (2\lambda+1)^d \text{ whenever } \lambda \in \mathbb{N}
		\label{natural}
	\end{equation}
	Now, let $\lambda > 1$ be any real number. Letting $\lfloor \lambda \rfloor$ denote the greatest integer less than or equal to $\lambda$, it is clear that \[N^{\infty}(\lambda) = N^{\infty}(\lfloor \lambda \rfloor) \leq N^{\infty}(\lfloor \lambda \rfloor + 1) = (2\lfloor\lambda\rfloor + 3)^d\]
	where the last equality follows from (\ref{natural}) with $\lfloor \lambda \rfloor + 1$ instead of $\lambda \in \mathbb{N}$.\\
	Finally, we have \[N^{\infty}(\lambda) \leq (2\lfloor\lambda\rfloor + 3)^d \leq (2\lambda + 3)^d \] where the second inequality follows from the fact that $\lfloor \lambda \rfloor \leq \lambda$ and $x \to (2x+3)^d$ is an increasing function on $\mathbb{R}_+$.
\end{proof}

With Lemma \ref{sup_norm} in hand, the proof of Proposition \ref{lattice_count} follows easily.

\noindent\textbf{Proof of Proposition \ref{lattice_count}:}
	As all norms are equivalent on $\mathbb{R}^d$, we have positive constants $m,M$ such that for any $x \in \mathbb{R}^d$
	\begin{equation} 
		mp(x) \leq p_{\infty}(x) \leq Mp(x) \label{valent_norms}
	\end{equation}
	We make a simple observation: for any $\lambda > 0$ 
	\begin{equation}
		B_p(\lambda) \subset B_{\infty}(M\lambda) \label{inclusion}
	\end{equation}
	This follows immediately from the second inequality in (\ref{valent_norms}).
	\\From (\ref{inclusion}) and the definitions of $N^p$ and $N^{\infty}$ we get that  for any $\lambda > 0$, \begin{equation}
		N^p(\lambda) \leq N^{\infty}(M\lambda)
		\label{link}
	\end{equation}
	Applying Lemma \ref{sup_norm} to $N^{\infty}(M\lambda)$, we get that for any $\lambda > 0$
	\begin{equation}
		N^{\infty}(M\lambda) \leq (2M\lambda + 3) ^d
		\label{ineq}
	\end{equation}
	Finally, (\ref{link}) and (\ref{ineq}) combine to give that for any $\lambda > 0$, \[N^p(\lambda) \leq N^{\infty}(M\lambda) \leq  (2M\lambda +3)^d\]
	And that is enough to show that  $N^p(\lambda) \in O(\lambda ^d)$ as $\lambda \to \infty$ \qed

\section{Acknowledgement} 
We are grateful to Antareep Saud for several useful discussions 
related to this project. 
The second author is also grateful to Vamsi Pingali for 
helpful discussions and comments and in particular pointing out the 
reference \cite{Getz_Ber} for the proof of Weyl's law.



\end{document}